\documentclass[12pt, a4paper]{amsart}
\usepackage{amsmath,amssymb,amscd,amsfonts}

\usepackage{ifthen}
\usepackage[T1]{fontenc}
\usepackage[utf8]{inputenc}
\usepackage[all]{xy}
\usepackage{graphicx}
\usepackage{enumerate}
\usepackage{xspace}
\usepackage{pdfsync}
\usepackage{epic}
\usepackage{dsfont}
\usepackage{tikz-cd}
\usepackage{hyperref}

\newtheorem{thm}{Theorem}[section]

\newtheorem{prop}[thm]{Proposition}
\newtheorem{cor}{Corollary}
\newtheorem{conj}{Conjecture}[section]
\newtheorem{ques}{Question}

\theoremstyle{definition}
\newtheorem{defn}{Definition}[section]
\newtheorem{exam}{Example}[section]

\theoremstyle{remark}
\newtheorem{rem}[thm]{Remark}

\DeclareMathOperator{\Bl}{Bl}
\DeclareMathOperator{\Hom}{Hom}

\begin{document}
\title[]{Around Homology Planes: Old and new}

\author{Rodolfo Aguilar Aguilar}
\address{\parbox{\linewidth}{University of Miami, Coral Gables, USA}}
%\thanks{The first author is partially supported by the Bulgarian Ministry of Education and Science, Scientific Programme "Enhancing the Research Capacity in Mathematical Sciences (PIKOM)", No. DO1-67/05.05.2022.}

\email{\href{mailto:aaguilar.rodolfo@gmail.com}{aaguilar.rodolfo@gmail.com}}
\urladdr{\url{https://sites.google.com/view/rodolfo-aguilar/}}

\date{\today}

\begin{abstract} 
We survey some old and new results about acyclic (affine) complex surfaces, also called homology planes. We ask several questions and leave open directions for future research.
\end{abstract}

\maketitle

\tableofcontents

\section{Introduction}
This note is to extend a talk given in Sunny Beach on August 2023 during the conference ``Generalized and Symplectic Geometry''. The talk was intended to give a survey around the classification of acyclic affine complex algebraic surfaces and to highlight some connections with topology and symplectic geometry.

Most of the note is expository although there are some new results. These new results are in the topology section and they will be explained in more detail below. 

In writing these notes, I use the opportunity to present some open questions that arose in discussions with people in the conference and outside. There are already some very good surveys about homology planes, see \cite{Z98} and \cite{M01}. For connections with $\mathbb{A}_1$-homotopy theory see also \cite{Asok21}. I try by no means to be exhaustive, and the contents here only reflects the interest of the author. He apologizes if he forgot to mention or give some credits.

The paper is organized as follows: in Section \ref{S:Context}, we present the cancellation problem as a motivation for studying homology planes and present some classification results.

In Section \ref{S:exotic}, we present some algebraic and analytic cancellation Theorems and show how this give exotic algebraic and analytic structures on $\mathbb{C}^n$ for $n\geq 3$. After this, we mention how we can obtain exotic symplectic structures on some affine real spaces. We end the section by presenting a result about exotic balls.

In Section \ref{S:MoreSy}, more questions and symplectic techniques are presented.

In Section \ref{S:Topology}, we study the relation of the fundamental group of some homology planes and the fundamental group of the $3$-manifold that can be thought as its boundary. We present a new result about the structure of this variety in Theorem \ref{thm:AN} and a partial computation of its character variety in Theorem \ref{thm:AS}.

Finally, in Section \ref{S:Alg}, we survey some results and questions in higher dimension and we finish presenting some open question about the cancellation problem.
 
\subsection*{Acknowledgements}
The author is grateful to F. Campana, J. Fernández de Bobadilla, L. Katzarkov,  M. McLean, A. Némethi, D. Pomerleano,  P. Popescu-Pampu, N. Saveliev, J. Svoboda, T. Syed and M. Zaidenberg for interesting discussions around the topics treated here. He also wants to thank the organizers and participants of the conference ``Generalized and Symplectic Geometry''.

He is thankful to A. Dubouloz, M. Golla, T. Pe{\l}ka, O. \c{S}avk,  for reading and commenting preeliminaries versions of this note. 

%Mention:

%That I am leaving some questions open or future directions. 

%Mention that I want to study further relation between the representations of the boundary and those in the interior.

\section{Context}\label{S:Context}

\subsection{Cancellation property}
Consider smooth affine varieties $X,Y,Z$.  We say that the triple\footnote{In this note, we will be mainly interested in the case when $Z=\mathbb{C}^n$ for some $n\geq 1$ and either $X$ or $Y$ will belong to some special class of quasi-projective varieties, for example, $Y$ is of log-general type. In this case, we don't fix the triple but only the special class of varieties, see Section \ref{S:exotic}.} $(X,Y,Z)$ satisfy the \emph{cancellation property} if, whenever we have an isomorphism $X\times Z \cong Y\times Z$ in the category of algebraic varieties, it induces a commutative diagram:

\[\begin{tikzcd}
X\arrow[d]\arrow[r,"\cong"]\times Z & Y \times Z \arrow[d] \\
X\arrow[r, "\cong"] & Y
\end{tikzcd}\]

As stated in \cite{R71}, it was believed around 1970, that every smooth contractible (affine), rational surface  was isomorphic to $\mathbb{C}^2$. 

This belief, together with the motivation of the cancellation problem, lead C.P. Ramanujam to find the following example.

\begin{exam}[\cite{R71}]\label{ex:Ram} Consider a plane conic $C_1$ and a cuspidal cubic $C_2$ intersecting as $C_1\cdot C_2=P +5Q$.  Let $\bar{X}=\Bl_P\mathbb{P}^2$ and denote the strict transforms of $C_1$ and $C_2$ by $C_1',C_2'$ respectively. Let $X=\bar{X}\setminus C_1'\cup C_2'$. C. P. Ramanujam showed that $X$ is affine contractible and $X\not \cong \mathbb{C}^2$.
\end{exam}

\begin{figure}[h]
\centering
\includegraphics[scale=.5]{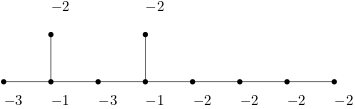}

\caption{Dual graph of the minimal embedded resolution of $C_1'\cup C_2'$}
\end{figure}

Moreover, using the methods for computing the fundamental group developed in \cite{M61}, he found a nice topological characterization of the complex affine plane. Before stating it, let us introduce a definition.

If $X$ is a smooth affine surface, we can consider a compactification $\bar{X}$ with $D:=\bar{X}\setminus X$ a SNC-divisor.  Take a closed regular tubular neighborhood $U$ of $D$ in $\bar{X}$.  The \emph{fundamental group at infinity} is $\pi_1(\partial U)$. 

\begin{thm}[\cite{R71}]\label{thm:Ram}
Let $X$ be a smooth affine surface which is contractible and simply connected at infinity. Then $X$ is isomorphic to $\mathbb{C}^2$ as an algebraic variety.
\end{thm} 
 
Ramanujam left open the following question. We will review the solution in Section \ref{ss:aES}. 
\begin{ques}[\cite{R71}] Is $\mathbb{C}^3\cong X\times \mathbb{C} $ as an algebraic variety?
\end{ques}

\subsection{Classification}

\begin{defn} Let $X$ be a smooth surface such that $H_i(X,\mathbb{Z})=0$ for $i>0$. We call $X$ a homology plane. If $H_i(X,\mathbb{Q})=0$ for $i>0$, we call $X$ a $\mathbb{Q}$-homology plane.
\end{defn}
\begin{exam} The Ramanujam's surface $X$ is a homology plane.
\end{exam}

Using the classification of smooth projective surfaces, the Miyaoka inequality \cite{Mi84} together with some technical computations,  Gurjar and Shastri proved the following:

\begin{thm}[\cite{GS89-1},\cite{GS89-2}] \label{thm:GS}
Homology planes are rational.
\end{thm}

The classification of homology planes follows the lines of the classification of surfaces via the (logarithmic) Kodaira dimension, that we introduce next.

\begin{defn}  Let $X$ be a smooth quasi-projective surface and let $(\bar{X},D)$ be a SNC-compactification. We define the \emph{logarithmic Kodaira dimension} $\bar{\kappa}(X)\in \{-\infty, 0,1,2\}$ as the unique value that satisfies that, for some positive constants $\alpha, \beta$ we have 
$$\alpha m ^{\bar{\kappa}(X)}\leq \dim H^0(\bar{X}, m(K_{\bar{X}}+D))\leq \beta m^{\bar{\kappa}(X)} $$
where $m$ is a sufficiently large and divisible positive integer.
\end{defn}

The invariant $\bar{\kappa}(X)$ does not depend on the chosen compactification of $X$. See \cite{I77}.

\begin{figure}[h]
       \centering
     \begin{tikzpicture}[scale=.8]
%\draw (.8,-1.5) node [below] {\scriptsize $q$};

\draw (-.4,2)  -- (1,-.5)node[right]{\tiny $L_n$};
\draw (-.2,2) -- (.5,-.5)node{\tiny $L_{n-1}$};
\draw (.2,2) -- (-.5,-.5)node[above,right]{\tiny $L_2$};
\draw (.4,2) -- (-1,-.5)node[above,right]{\tiny $L_1$};

\draw (-1,0) -- (1,0) node[right]{\tiny $L_{n+1}$};

%\foreach \Point in {(.4,0),(-.4,0),(0,.4),(0,-.4)}{
%    \node at \Point {$\cdot$};
%}

%\draw [red] (-.4,.6) to[out=-150,in=100] 
%(-.5,.3) [arc arrow=to pos 1 with length .5mm] %;

%\draw [red] (-.7,0) to[out=-150,in=100] 
%(-.8,-.3) [arc arrow=to pos 1 with length .5mm] ;
\end{tikzpicture}
\caption{$L(1,5)$ in the projective plane}
\label{fig:LinesL1}
    
\end{figure}
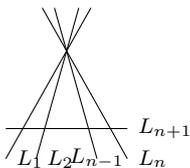

When $\bar{\kappa}(X)\leq 1$, we have an easy description for the homology plane $X$. Recall that a quasi-pencil $L(1,n+1)$ is a set of lines in $\mathbb{P}^2$ having $n$-lines passing through a point and another line in general position. See Figure \ref{fig:LinesL1}.

\begin{thm}[\cite{F82}, \cite{GM88}, \cite{tDP90}] \label{thm:ClassLess1}
Let $X$ be a homology plane.
\begin{enumerate}
\item If $\bar{\kappa}(X)=-\infty$ then $X\cong \mathbb{C}^2$.
\item There are no homology planes with $\bar{\kappa}(X)=0$.
\item All homology planes with $\bar{\kappa}(X)=1$ arise from a quasi-pencil of lines $L(1,n)$. 
\end{enumerate}
\end{thm}

For a homology plane $X$ with $\bar{\kappa}(X)=1$, its construction depends of certain integer pairs $(a_1,b_1),\ldots, (a_n,b_n)$ satisfying 
\begin{equation}\label{eq:determinant}
a_1\cdots a_n +\sum_{i=1}^n a_1\cdots b_i \cdots a_n =\pm 1.
\end{equation}
It is obtained from $L(1,n)$ by blowing-up over some singular points and removing some $(-1)$-curves. For its precise construction, see \cite{M01} or \cite[Section 2.2]{AS22}.

%Its dual graph looks as in  Figure \ref{fig:DualGraph}. We will study it more in detail in section \ref{S:Topology}.
%\begin{figure}[h]
%\centering
%\includegraphics[scale=.25]{dgs.png}

%\caption{Dual graph of homology planes with $\bar{\kappa}(X)=1$.}
%\label{fig:DualGraph}
%\end{figure}

Now, when $\bar{\kappa}(X)=2$ we say that $X$ is of log-general type. Amazingly, due to the constrains imposed by $X$ being a homology plane, there is a complete (conjectural) classification. 

Let us begin first with some conjectures.

\begin{conj}\label{conj:Koras} Let $X$ be a homology plane of log general type and let $(\bar{X},D)$ be a minimal SNC completion of $X$. Then 
\begin{enumerate}
\item (tom Dieck-Petrie conjecture) Every homology plane arises from finitely many configurations of lines and conics in $\mathbb{P}^2$.
\item (Koras conjecture) The automorphism group of $X$ has at most six elements, 
\item (Flenner-Zaidenberg Rigidity conjecture) $H^i(\mathcal{T}_{\bar{X}}(-\log D))=0$ for all $i\geq 0$, 
\item (Palka negativity conjecture) $\kappa(K_{\bar{X}}+\frac{1}{2}D)=-\infty$.
\end{enumerate}
\end{conj}
 
In \cite{P23}, T. Pe\l ka completely classified the $\mathbb{Q}$-homology planes assuming Palka negativity conjecture. He uses the theory of almost minimal model of Palka \cite{P19}. Indeed, assuming that $\kappa(K_{\bar{X}}+\frac{1}{2}D)=-\infty$, one obtains that a minimal model $(\bar{X}_{\text{min}},\frac{1}{2} D_{\text{min}})$ is a log-Mori fiber space over either a curve or a point. The first case was studied by \cite{MS91}, although Pe\l ka gives a self-contained proof. The second case is the main part of his work. One can then state his main theorem (we will state it only for $\mathbb{Z}$-homology planes).

\begin{thm}[\cite{P23}]
Assume (4) in Conjecture \ref{conj:Koras} (Palka negativity conjecture), then it implies (1),(2) and (3). In particular, the $\mathbb{Z}$-homology planes of log-general type can be arranged in finitely many discrete series, each obtained in a uniform way from an arrangement of lines and conics on $\mathbb{P}^2$.
\end{thm}

\section{Exotic structures}\label{S:exotic}
\subsection{Algebraic}\label{ss:aES}
Using the theory developed in \cite{I77}, Iitaka-Fujita proved the following cancellation theorem.
\begin{thm}[\cite{IF77}]\label{thm:IF}
If we have an algebraic isomorphism $X\times \mathbb{C}^n\cong Y\times \mathbb{C}^n$ and $\bar{\kappa}(Y)\geq 0$ then the cancellation property holds: $X\cong Y$.
\end{thm}

By taking $Y$ to be the Ramanujam surface of Example \ref{ex:Ram} (it can be shown that $\bar{\kappa}(Y)=2$), we can obtain algebraic exotic structures on $\mathbb{C}^n$ for $n\geq 3$.
\begin{cor}\label{cor:AlgExotic} The space $\mathbb{C}^n$ admits an exotic algebraic structure for $n\geq 3$.
\end{cor}

Indeed, using the $h$-cobordism Theorem one can show that $Y\times \mathbb{C}^n$ is diffeomorphic to $\mathbb{R}^{2n+4}$. By taking $X=\mathbb{C}^2$ and if we suppose that $X\times \mathbb{C}^n\cong Y\times \mathbb{C}^n$,  we use the cancellation Theorem \ref{thm:IF} of Iitaka-Fujita, which would give the contradiction that the Ramanujam surface is isomorphic to $\mathbb{C}^2$.

\subsection{Analytic}
Using properties of non-singular quasi-projective varieties of log-general type (e.g. every dominant holomorphic map $X\to Y$ from a quasi-projective manifold to a manifold of log-general type is regular), Zaidenberg proved the following \emph{analytic} cancellation theorem.
\begin{thm}[\cite{Z93}]\label{thm:ZaiCan} Suppose that $X$ and $Y$ are smooth surfaces of log-general type, this is $\bar{\kappa}(X)=\bar{\kappa}(Y)=2$, then $X,Y,\mathbb{C}^n$ satisfy an \emph{analytic} cancellation property.
\end{thm} 

Moreover, he constructed two infinite families of contractible homology planes of log-general type pairwise non-isomorphic. One family generalized the Ramanujam surface.

Using these infinite families, his analytic cancellation theorem \ref{thm:ZaiCan} and the same argument in the proof of Corollary \ref{cor:AlgExotic}, he obtained the following Corollary.

\begin{cor}[\cite{Z93}] There is a countable infinite set of analytic exotic structures on $\mathbb{C}^n$ for $n\geq 3$. 
\end{cor}

Moreover, in \cite{FZ92}, by using that homology planes of log-Kodaira dimension one can have moduli, they obtained an \emph{uncountable} number of exotic structures on $\mathbb{C}^n$ for $n\geq 3$.

\subsection{Symplectic}
In \cite{SS05}, the authors studied the symplectic topology of the Ramanujam surface. More precisely, they obtained the following.
\begin{thm}[\cite{SS05}]\label{thm:SS05} Let $X$ be the Ramanujam's surface, the product $X^m\cong \mathbb{R}^{4m}$ induces an exotic symplectic structure on $\mathbb{R}^{4m}$ which is convex at infinity.
\end{thm}

As explained in \cite{S08}, this can be shown by showing that the symplectic cohomology of $X^m$ does not vanish. 

To construct an infinite countable number of exotic symplectic structures on $\mathbb{R}^{2m}$ for $m>3$ McLean started from a different example. He used the so-called Kaliman modification \cite{K94} to construct Lefschetz fibrations. As the cancellation theorem are no longer valid in the symplectic category, he used the number of idempotents in symplectic homology as invariant to distinguish between different exotic structures.

\begin{thm}[\cite{M09}] Let $m>3$. Then for each $\mathbb{R}^{2m}$ there are infinitely many finite type Stein manifolds which are pairwise distinct as symplectic manifolds.
\end{thm}

\begin{ques} Is it possible to construct uncountably many exotic symplectic structures on $\mathbb{R}^{2m}$ for $m>3$?
\end{ques}

Note that one needs a stronger invariant than symplectic cohomology if one wants to use the deformation examples of Flenner-Zaidenberg.

\subsection{Contractible four-manifolds}

An important problem in low-dimensional topology is to determine which homology $3$-spheres bound contractible $4$-manifolds. See \cite{Sa23}.

A \emph{Mazur} manifold, is a contractible, compact, smooth $4$-dimensional manifold with boundary, which is not diffeomorphic to the standard $4$-ball and such that in its handle decomposition it contains exactly three handles: one $0$-handle, one $1$-handle and one $2$-handle.

Let us review briefly the construction of Mazur and Poénaru manifolds as in \cite[Sec. 2.2]{Sa22} which is closer to the construction by V. Poénaru \cite{P60}. Consider a disk $D$ smoothly embedded in $B^4$ and take a tubular neighborhood $U$ of $D$. In the boundary $\partial U$, which is diffeomorphic to $S^1\times S^2$, make an integral surgery in a knot $J$ with fusion number $0$ such that the $3$-manifold obtained is a homology sphere. The resulting $3$-manifold bounds a Mazur manifold, see \cite[Lemma 2.1]{Sa22}.

In this construction, if we replace the disk $D$ by a ribbon disk $D\subset B^4$ and the knot $J$ by a ribbon knot with fusion number $n\geq 1$, the resulting $3$-manifold bounds a contractible four manifold with one $0$-handle, $n+1$ $1$-handles and $n+1$ $2$-handles. This is called a \emph{Poénaru manifold}.

In order to obtain this ribbon knot, one fixes a compactification $(\bar{X},D)$ of a homology plane $X=\bar{X}\setminus D$ and consider the dual graph of $D$. This can also be seen as the surgery diagram of the boundary of a tubular neighborhood of $D$ as follows: to every vertex associate an unknot, link two unknots if and only if there is an edge connecting the corresponding vertices and decorate the unknots with the weights of the corresponding vertices. If by applying Kirby calculus one ends with a ribbon knot, one uses {\c{S}}avk's construction of Poénaru manifolds. 

By using the vast quantity of example in the literature of homology planes, in particular some examples of Zaidenberg generalizing the Ramanujam surface and contractible homology planes of log-Kodaira dimension one and using Kirby calculus, one can show the following Theorem.

\begin{thm}[\cite{AS22}]\label{thm:AgSav} There are four infinite families of homology $3$-spheres that bound both homology planes and Mazur/Poénaru manifolds.
\end{thm}

In \cite{AS22}, these homology $3$-spheres bounding both homology planes and Mazur/Poénaru manifolds  are called Kirby-Ra\-ma\-nu\-jam spheres. In Theorem \ref{thm:AgSav}, one infinite family consists of Mazur manifolds and the other three of Poénaru manifolds.

\section{More symplectic properties}\label{S:MoreSy}

\subsection{Search for more exotic symplectic structures}

One can ask if by using other examples of contractible homology planes one can obtain more exotic symplectic structures on $\mathbb{R}^{4m}$. This is the content of the unpublished work of D. Jackson-Hanen \cite{JH14}. 

Moreover, partially motivated by this kind of quest, Ganatra-Pomerleano showed the following Theorem.

\begin{thm}[\cite{GP20}]\label{thm:GP20} Let $X$ be a smooth affine variety and let $(\bar{X},D)$ be a SNC-compactification. Then there is a multiplicative spectral sequence converging to the symplectic cohomology ring 
$$\{E_r^{p,q},d_r\}\to SH^*(X)$$
whose first page is isomorphic as rings to the logarithmic cohomology of $(\bar{X},D)$.
\end{thm} 

For the definition of logarithmic cohomology of a pair $(\bar{X},D)$, see \cite{GP21}.

\begin{ques}

\begin{enumerate}
\item Can one use this spectral sequence to recover the Seidel-Smith result? 
\item Can one use the classification of Pe\l ka to obtain more exotic structures?

\end{enumerate}
\end{ques}

\subsection{A question for the contact structure of the boundary} 
If $X=\bar{X}\setminus D$ is a homology plane and we let $U$ be a closed regular tubular neighborhood of $D$, then the boundary of $\bar{X}\setminus U$ inherits a contact structure.

We can say that an oriented contact $3$-manifold $(M,\xi)$ is Ramanujam fillable if it is contactomorphic to the contact boundary of a homology plane. 

We can try to extend the following result of singularity theory.

\begin{defn}Let $(M,\xi)$ be an oriented contact manifold. We say it is \textit{Milnor fillable} if it is contactomorphic to the contact boundary of an isolated complex-analytic singularity $(\mathcal{X},x)$.
\end{defn}

\begin{thm}[\cite{CNPP06}]\label{thm:CNPP} Let $(M,\xi)$ be any three dimensional oriented manifold. Then it admits at most one Milnor fillable structure up to contactomorphism.
\end{thm}
 
\begin{ques}
Does an analog of the Caubel, N\' emethi, Popescu-Pampu theorem \ref{thm:CNPP} holds in this setting (probably restricting the Kodaira dimension)?
\end{ques}

In fact, as shown in \cite[Cor 1.10]{P23}, there exists examples of non (algebraically) isomorphic but diffeomorphic homology planes $X=\bar{X}\setminus D$ with the same dual graph for $D$. These are sometimes called Zariski-pairs, see \cite{ABCT08}. This indicates that additional information to the dual graph is needed, we ask if the contact structure is enough for this.  Moreover, similar questions and techniques have been studied in the symplectic category for rational cuspidal curves. See \cite{GS22, GK23}.

\subsection{Invariance of log-Kodaira dimension}
Let us gather here more questions related to homology planes. 

In the compact setting, the invariance of the Kodaira dimension under diffeomorphism is shown by Friedman-Morgan \cite{FM94}. In the symplectic category we have the following result of McLean.

\begin{thm}[\cite{M14}] The logarithmic Kodaira dimension is a symplectic invariant of homology planes.
\end{thm} 

To prove this, he uses Theorem  \ref{thm:Ram}, \ref{thm:ClassLess1} and the theory he develops about the symplectic invariance of uniruledness.

\begin{ques} (McLean): Is the logarithmic Kodaira dimension a topological invariant of homology planes?
\end{ques}
\subsection{Mirror symmetry}

%\textit{Problem}: Compute the derived category of a homology plane $X$. 
When $X$ is a contractible homology plane of log-Kodaira dimension one, it is shown in \cite{tDP90} that they can arise as hypersurfaces in $\mathbb{C}^3$. It could be interesting to compute their mirror using the techniques in \cite{AAK16}.

Moreover, as they are contractible, one can wonder what kind of invariants could be non-trivial in the mirror. Recall that for example, in the Ramanujam surface, symplectic cohomology is non-trivial. 

The following question is due to D. Pomerleano. 
\begin{ques} Is it possible to use the homology planes to construct exotic $(\mathbb{C}^{*})^r$ structures? 
\end{ques}

Note that, if one wants to follow the proof of Corollary \ref{cor:AlgExotic}, one would need to use the s-cobordism Theorem, instead of the h-cobordism.

%[Pomerleano] connection with toric varieties, 

%\begin{prop}[\cite{AD08}] Let $X$ be a contractible homology plane. Then every vector bundle on $X$ is isomorphic to a trivial bundle.
%\end{prop}

\section{Topology}\label{S:Topology}
\subsection{Fundamental groups}

Motivated by the question, often attributed to Serre, of what groups can appear as fundamental groups of complex projective manifolds, we can ask:
\begin{enumerate}
\item What groups can appear as fundamental groups of homology planes,  
\item What groups can appear as fundamental groups at infinity for homology planes, 
\item What is the relation between these two?  
\end{enumerate}

A first constrain in given by the following theorem of S. Orevkov.

\begin{thm}[\cite{O97}]\label{thm:Orev} Let $X=\bar{X}\setminus D$ be a homology plane and for $D\subset U$ a tubular neighborhood as above then the three-manifold $\partial U$ is not a Seifert fibered homology sphere.
\end{thm}

 An \emph{orbicurve group} is a group admitting a presentation of the following type:
\begin{align*}
&\langle a_1,b_1,\ldots,a_g,b_g,c_1,\ldots,c_s\mid \\
& \quad \quad c_1^{r_1}=\ldots=c_s^{r_s}=1, \prod_{i=1}^g [a_i,b_i]c_1\ldots c_r=1\rangle 
\end{align*}
\begin{thm}[\cite{AA21}] \label{thm:AA21} There exists infinitely many homology planes of log-general type with infinite fundamental group. All of them admit a quotient to an infinite orbicurve group.
\end{thm}

For the proof, a presentation of the fundamental group of the homology plane is obtained. Then, some simpler quotient (an orbicuve group) which is known to present an infinite group is deduced.

\begin{ques}\label{ques:fund} Does every infinite fundamental group of a homology plane admits a surjective morphism to an infinite orbicurve group?
\end{ques}

By a Theorem of Arapura \cite{A97}, extended to orbicurves in  \cite{ABCAM13}, the morphisms from quasi-projective varieties to  orbicurves are encoded in the irreducible components of the so-called \textit{cohomology jump loci}. For a quasi-projective variety $X$, this loci consists of elements $\rho$ in $H^1(X,\mathbb{C}^\ast)$ such that, the associated local system $V_\rho$ has non-vanishing first cohomology $H^1(X, V_\rho)\not = 0$.

 Now, from the construction of homology planes coming from an arrangement of lines $\mathcal{A}$, every point of multiplicity three of higher will give a map to $\mathbb{P}^1$ with some orbifold structure. However, to find the examples in Theorem \ref{thm:AA21}, we used the complete description of the cohomology jump loci as given by Suciu in \cite{S01}. In particular, we used maps not coming from points of multiplicity three or higher, but others reflecting the global structure of $\mathcal{A}$ plus some carefully chosen numerical data that depend on the homology plane.

We do not know if the question \ref{ques:fund} has a positive answer even if we assume the negativity conjecture \ref{conj:Koras}, but for homology planes of Kodaira dimension one it has indeed a positive answer. See below remark \ref{rem:hpfund}.

If one would like to study the cohomology jump loci in the homology plane directly and not in the complement of an arrangement $\mathcal{A}$ in $\mathbb{P}^2$ a generalization of Arapura's result is needed for non-abelian representations, for example in $SL_2(\mathbb{C})$. Even simpler, we can try to study the $SU(2)$ character varieties. This is not easy, even for the representations of the $3$-homology sphere bounding a homology plane $X$, as we will see in subsection \ref{ss:char}. For more of representations of quasi-projective groups on $SL_2(\mathbb{C})$, see \cite{CS08, PS20}.

When $X$ is a homology plane of log-Kodaira dimension one. It is easy to compute its fundamental group. 

\begin{prop} Let $X$ be a homology plane of logarithmic Kodaira dimension one coming from $L(1,n+1)$ with data $(a_1,b_1),\ldots, (a_n,b_n)$ then a presentation for its fundamental group is given by 
\begin{equation}\label{eq:PresHomPla}
\langle \gamma_0,\ldots \gamma_{n}\mid \gamma_1\cdots \gamma_{n}=\gamma_0^{-1} , [\gamma_0,\gamma_i]=1, \gamma_i^{a_i}\gamma_0^{b_i}=1 \text{ for }i=1,\ldots,n\rangle
\end{equation}
\end{prop}
\begin{rem}\label{rmk:Determinant} Note that from (\ref{eq:PresHomPla}), we obtain that $H_1(X)$ is trivial if and only if 
$$\det \left(\begin{array}{rrrr}
a_1-b_1 & -b_1 & \ldots & -b_1\\
-b_2 & a_2-b_2 & & -b_2 \\
\vdots &  & \ddots & \\ 
-b_n & -b_n & & a_n-b_n

\end{array} \right)=\pm 1 $$
which is precisely the condition in (\ref{eq:determinant}).
\end{rem}

\begin{rem}\label{rem:hpfund}  Note that all homology planes obtained from $L(1,3)$ are contractible. Indeed, note that from (\ref{eq:PresHomPla}), more specifically, from $[\gamma_0,\gamma_2]=1$ and $\gamma_0^{-1}=\gamma_1\gamma_2$ we obtain that $[\gamma_1,\gamma_2]=1$ and hence $\pi_1(X)$ is abelian. By the condition (\ref{eq:determinant}) and by Remark (\ref{rmk:Determinant}), we obtain that $\pi_1(X)$ is trivial. 

Moreover, when $n\geq 3$, the homology planes obtained from $L(1,n+1)$ by the above procedure are never contractible. Indeed, it admits an orbifold morphism to $\mathbb{P}^1$ with $n$ marked points with weights $a_1, \ldots, a_n$. 

%is reduces to 
%$$\langle \gamma_1,\gamma_2 \mid \gamma_1^{a_1}\gamma_2^{b_2} \rangle $$
\end{rem}
\begin{rem} Note that the presentation in (\ref{eq:PresHomPla}) is the presentation of the fundamental group of a Seifert fibered homology sphere.
\end{rem}

\subsection{Boundary manifolds}

\begin{figure}[h]
\centering
\includegraphics[scale=.25]{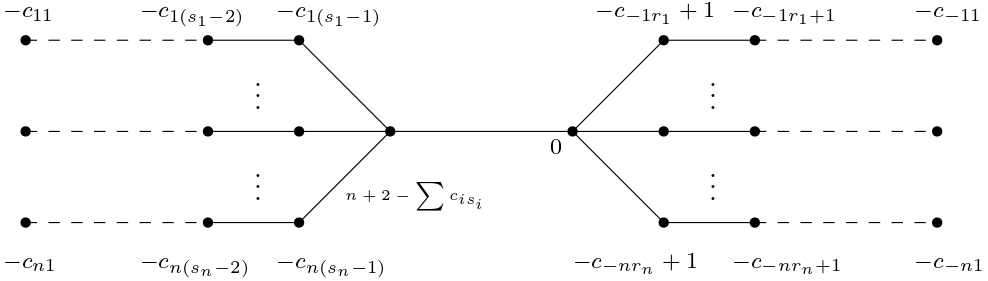}

\caption{Splice of the Seifert spheres $\Sigma(a_1,\ldots,a_n, a_1\cdots a_n \pm 1)$ along its last fiber.}
\label{fig:DualGraph}
\end{figure}

When $X=\bar{X}\setminus D$ is a homology plane of log-Kodaira dimension one, the dual graph of $D$ is as in Figure \ref{fig:DualGraph} or it can have an extra linear segment of $-2$-curves corresponding to a series of blow-ups starting in a smooth point in the line $L_{n+1}$.

Recall Theorem \ref{thm:AgSav}. In \cite{AS22}, for an infinite family $(\bar{X}_n,D_n)$ of homology planes of logarithmic Kodaira dimension one, it was shown that if we let $\partial U_n$ be the boundary of a regular tubular neighborhood of $D_n$, then it is the splice of two Seifert fibered spheres $\Sigma(a_{1n},a_{2n},a_{1n}a_{2n}\pm 1)$ glued along the last fiber. That is, in each Seifert fibered homology sphere $\Sigma(a_{1n},a_{2n},a_{1n}a_{2n}\pm 1)$ a small neighborhood of the third twisted fiber is removed, and the boundaries are glued longitude-to-meridian.

Using the techniques in \cite{NR78}. Together with A. N\'emethi we obtained the following Theorem. 
\begin{thm}[A-N\'emethi '22]\label{thm:AN} Let $X$ be a homology plane with $\bar{\kappa}(X)=1$ and suppose that the dual graph of $D$ is as in Figure \ref{fig:DualGraph}. Let $\partial U$ be as above. Then the homology sphe $\partial U$ is the splice of $$\Sigma(a_1,\ldots,a_n, a_1\cdots a_n\pm 1)$$ along the last fiber. 
\end{thm}

It is easy to compute the fundamental group of $\partial U$ by using the methods of \cite{M61}. Moreover, by fixing explicit geometric meridians around the lines in $L(1,n)$, it is possible to compute the presentation in (\ref{eq:PresHomPla}) as a quotient from a presentation of $\pi_1(\partial U)$. We will see in the next section that it is however, very hard to compute and find relations between invariants of these four and three (real) manifolds.

\subsection{Character varieties}\label{ss:char}
Let $X$ be a compact manifold, the $SU(2)$ \emph{character variety} of $X$ is defined to be the topological space 
$$R(X)=\Hom(\pi_1(X),SU(2))/conj $$

For homology planes of log-Kodaira dimension, by Theorem \ref{thm:AN}, we have that the $3$-manifold $\partial U$ the splice of two Seifert spheres $$\Sigma(a_1,\ldots, a_n, a_1\cdots a_n \pm 1)$$ along the last fiber.

When the splice components have only three singular fibers, the character variety is studied by N. Saveliev in \cite{S98}. Interestingly, it turns out that it is simply the disjoint union of the two splice component. 

\begin{thm}[\cite{S98}]\label{thm:Sav98} Let $\Sigma$ be the splice of two Seifert spheres $\Sigma(a_1, a_2, a_1 a_2\pm 1)$ along its last fiber. Then $$R(\Sigma)=R(\Sigma(a_1, a_2, a_1 a_2+ 1))\coprod R(\Sigma(a_1, a_2, a_1 a_2 - 1))  $$
\end{thm}

Using this, he can conclude that instanton Floer homology groups $I_j(\Sigma)$ vanish for odd $j$.

Now, if we want to generalize this computations to the case when $n>2$, this is, if the splice components have more than $3$ singular fibers, the computations are much more complicated and only partial results can be obtained with similar computations. 

Denote $\Sigma_\pm= \Sigma(a_1,\cdots, a_n, a_1\cdots a_n \pm 1)$ and let $$\mathcal{R}(\Sigma_\pm)=\mathcal{R}_\pm \coprod \mathcal{R}(a_1,\ldots,a_n)$$ where $\mathcal{R}_+$ is the union of the components of $\mathcal{R}(\Sigma_+)$ consisting of those representations $\rho$ that send the generator corresponding to the fiber last fiber $a_1\ldots a_n +1$ to a non-trivial element. A similar description holds for $\mathcal{R}_-$.

The following description was obtained in discussions with N. Saveliev. The proof is analogous to that of Theorem \ref{thm:Sav98}.

\begin{thm}[A-Saveliev '23]\label{thm:AS} Let $\Sigma$ be the splicing of $\Sigma_+$ and $\Sigma
_-$ along the $n+1$ singular fiber. Then $$\mathcal{R}(\Sigma)\supset  E \coprod \left(\mathcal{R}_+ \coprod \mathcal{R}_-\right),$$
where $E$ is an $SO(3)$-fiber bundle over $\mathcal{R}(a_1,\ldots,a_n) \times \mathcal{R}(a_1,\ldots,a_n)$.\end{thm}

The difficult part comes from those representations of $\pi_1(\Sigma)$ that restrict to a central representation in one component and are irreducible in the other. Of course, this also complicates the computation of instanton Floer homology groups. 

It may happen that other invariants could give more insight in the relation between the geometric and algebraic structure. For example if the intersection matrix of $D$ is negative definite, the analytic lattice cohomology of A. N\'emethi and collaborators gives nice relations, see \cite{AN21}. However, if the matrix is no longer negative definite, it is difficult to conclude. It could be interesting if the invariants of \cite{GPPV20} could be computed (if defined) for $\partial U$.

\section{Algebra}\label{S:Alg}

\subsection{Higher dimensions}
It is a result from Fujita \cite{F82}, that homology planes are affine surfaces. In particular they are Stein. If we consider contractible algebraic varieties of higher dimension, Winkelmann showed that this is no longer the case.
\begin{thm}[\cite{W90}] There exists a contractible quasi-affine variety $X$ such that $X$ is non-Stein.
\end{thm}

The example is obtained by the quotient of a free action of the additive group $\mathbb{G}_a$ on $\mathbb{A}^5$. One may wonder if more examples of this type could be constructed.

\begin{thm}[\cite{ADF17}] Let $n\geq 3$. There exists contractible varieties $X_{2n}$ of dimension $2n$ which are not a quotient of affine space by a free action of a unipotent group. 
\end{thm}

The examples are given by removing a closed subscheme $E$ isomorphic to $\mathbb{A}^n$ from certain affine quadric $Q_{2n}$.

One may also ask of the higher dimensional generalization of the Theorem \ref{thm:GS} of Gurjar-Shasti  holds. This being probably too hard, can be weakened to the following question of F. Campana.

\begin{ques}[Campana] Let $X$ be a smooth, contractible affine variety $X$. Consider a smooth compactification $\bar{X}$ of $X$, is then $\bar{X}$ rationally connected? 
\end{ques}

Similar questions about compactifications of smooth contractible threefolds have been treated by T. Kishimoto and collaborators. See \cite{K04,K05}. %[Compactifications of contractible affine 3-folds into smooth Fano 3-folds with B2=2, On the compactifications of contractible affine
%threefolds and the Zariski Cancellation Problem]

%
%{Cancellation problem}

% The homotopy type conjecture, (weak form of the geometric $P=W$ conjecture) predicts that the dual boundary complex $\mathbb{D}\partial M_{\mathcal{B}}$ is homotopy equivalent to a sphere of dimension $d-1$.

The following comment was made by A. Dubouloz.

It is remarkable the existence of exotic algebraic structures on $\mathbb{C}^3$ that are rational, affine, of logarithmic Kodaira dimension $-\infty$ and such that it is not a product of a contractible surface with $\mathbb{C}$.

Probably, the most famous example is Russell's cubic $X$ defined by the equation ${x^2y+z^2+t^3+x=0}$ in $\mathbb{C}^4$. Topologically, it is a variety diffeomorphic to $\mathbb{R}^6$, but algebraically different from $\mathbb{C}^3$. The following open questions seem far from being answered.

\begin{ques} \begin{enumerate}
\item Is $X$ analytically isomorphic to $\mathbb{C}^3$?
\item Is $X\times \mathbb{C}$ algebraically/analytically isomorphic to $\mathbb{C}^4$?

\end{enumerate}
\end{ques}

In the symplectic context, it is shown in \cite{CM19} that $X$ is Stein deformation-equivalent to $\mathbb{C}^3$.
% mais il y a par contre de jolis résultats concernant cette variété (et des familles de variétés similaires) dans le contexte symplectique (cf. par exemple le papier de Casals et Murphy "Legendrian Fronts for Affine Varieties") qui peuvent également être mentionnés en complément de tous les aspects symplectiques en dimension complexe 2 que tu décris dans ton survol. 

\subsection{Cancellation problem}
To finish this note, let us conclude with an unexpected  connection , at least to the author, between the cancellation property and moduli spaces of representations on a punctured Riemann surface. 

Let $M_\mathcal{B}$ be the Betti moduli space of stable filtered regular $G$-local systems over a punctured Riemann surface. Denote the complex dimension of $M_\mathcal{B}$ by $d$. When $M_{\mathcal{B}}$ is very generic (see Section 1.2  in \cite{S23}), we have a strong form of Mellit's cell decomposition: 

\begin{thm}[\cite{S23}] $\mathcal{M}_B$ is decomposed into locally closed subvarieties of the form $(\mathbb{C}^*)^{d-2b}\times \mathcal{A}$, where $\mathcal{A}$ is stably isomorphic to $\mathbb{C}^b$.
\end{thm}

\textit{Expectation (Su)}: $\mathcal{A}$ should give in general a counter-example to the Zariski cancellation problem for dimension bigger than $3$ in characteristic zero.

\bibliography{HP}
\bibliographystyle{amsalpha}

\end{document}